\documentclass[11pt,a4paper,oneside,final]{article}

\usepackage{amssymb}
\usepackage{amsmath}
\usepackage{amsthm}
\theoremstyle{plain}
\newtheorem{theorem}{Theorem}
\newtheorem{corollary}[theorem]{Corollary}
\newtheorem{proposition}[theorem]{Proposition}
\newtheorem{lemma}[theorem]{Lemma}

\theoremstyle{definition}
\newtheorem{definition}[theorem]{Definition}
\newtheorem{example}[theorem]{Example}

\theoremstyle{remark}
\newtheorem{remark}[theorem]{Remark}

\begin{document}

\title{Lipschitzian Regularity of the Minimizing Trajectories
       for Nonlinear Optimal Control Problems\thanks{This 
       research was partially presented, as an oral communication, 
       at the international conference EQUADIFF 10, Prague, August 27--31, 2001.
       Accepted for publication in the journal Mathematics of Control, Signals, 
       and Systems (MCSS).}}

%\titlerunning{Lipschitzian regularity of the minimizing trajectories}

\author{Delfim F. M. Torres}

%\authorrunning{D. F. M. Torres}

\date{Department of Mathematics\\
         University of Aveiro\\
         3810-193 Aveiro, Portugal\\
         \texttt{ fax: +351 234382014\\
                  e-mail: delfim@mat.ua.pt
                }
        }
           
% http://www.mat.ua.pt/delfim
           
\maketitle
           
\begin{abstract}
We consider the Lagrange problem of optimal control
with unrestricted controls and
address the question: under what conditions we can assure
optimal controls are bounded? This question is related to the one
of Lipschitzian regularity of optimal trajectories, and the answer
to it is crucial for closing the gap between the conditions
arising in the existence theory and necessary optimality
conditions. Rewriting the Lagrange problem in a parametric form, we
obtain a relation between the applicability conditions 
of the Pontryagin maximum principle to the later problem
and the Lipschitzian regularity conditions for the original problem.
Under the standard hypotheses of coercivity of the existence theory, 
the conditions imply
that the optimal controls are essentially
bounded, assuring the applicability of the classical 
necessary optimality conditions like the
Pontryagin maximum principle. The result extends previous
Lipschitzian regularity results to cover optimal control problems with
general nonlinear dynamics.
\end{abstract}


%%%%%%%%%%%%%%%%%%%%%%%%%%%%%%%%%%%%%%%%%%%%%%%%%

\section{Introduction}

Given a Lagrangian $L : \mathbb{R} \times
\mathbb{R}^n \times \mathbb{R}^r \longrightarrow \mathbb{R}$, a
dynamical equation $\dot{x}(t) =
\varphi\left(t,\,x(t),\,u(t)\right)$, and boundary conditions $x(a) =
A$, $x(b) = B \in \mathbb{R}^n$, we consider the problem
of finding a control $u(\cdot) \in
L_1\left([a,\,b];\,\mathbb{R}^r\right)$ such that the
corresponding state trajectory $x(\cdot) \in
W_{1,\,1}\left([a,\,b];\,\mathbb{R}^n\right)$ of the dynamical
equation satisfies the boundary conditions, and the pair
$\left(x(\cdot),\,u(\cdot)\right)$ minimizes the functional
$J[x(\cdot),\,u(\cdot)] := \int_a^b L\left(t,\,x(t),\,u(t)\right) 
\mathrm{d}t$.
We establish Lipschitzian regularity conditions for the minimizing
trajectories of such optimal control problems. Lipschitzian regularity has a
number of important implications. For example in control
engineering applications, where optimal strategies are implemented
by computer, the choice of discretization  and numerical
procedures depends on minimizer regularity \cite{MR86k:49002,MR2002e:49058}.
Lipschitzian regularity of optimal trajectories also precludes occurrence of
the undesirable Lavrentiev phenomenon
\cite{MR2000f:49011,MR97g:49014,mania,JFM53.0481.02} and provides the validity of
known necessary optimality conditions under hypotheses of
existence theory \cite{MR2001h:49003}. 
The techniques of the existence theory use
compactness arguments which require to work with measurable
control functions from $L_p,\,1 \le p < \infty$
\cite{MR85c:49001}. On the other hand, standard necessary
conditions for optimality, such as the classical Pontryagin maximum
principle \cite{MR29:3316b}, put certain restrictions 
on the optimal controls --
namely, \emph{a priori} assumption that they are essentially
bounded. Examples are known, even for polynomial
Lagrangians and linear dynamics \cite{MR86k:49002}, for which
optimal controls predicted by the existence theory are unbounded
and fail to satisfy the Pontryagin maximum principle
\cite{MR85m:49051}. If we are able to assure that a minimizer
$(\tilde{x}(t),\,\tilde{u}(t))$, $ a \le t \le b$, of our problem
is such that $\tilde{u}(\cdot)$ is essentially bounded, then the
solutions can be identified via the Pontryagin maximum principle.
As far as $\varphi(t,\,\tilde{x}(t),\,\tilde{u}(t))$ is bounded, it
also follows that the optimal trajectory $\tilde{x}(\cdot)$ is
Lipschitzian. Similarly, the Hamiltonian adjoint multipliers
$\tilde{\psi}(\cdot)$ of the Pontryagin maximum principle turn out to
be Lipschitzian either. Thus,  regularity theory justifies
searching for minimizers among extremals and establishes a weaker
form of the maximum principle in which the Hamiltonian adjoint 
multipliers are not required to be absolutely continuous but
merely Lipschitzian.

The study of Lipschitzian regularity conditions has received few
attention when compared with existence theory or necessary
conditions, which have been well studied since the fifties and
sixties. The question of Lipschitzian regularity, for the general
Lagrange problem of optimal control, seems difficult, and
attention have been on particular dynamics. Most part of results
in this direction refers to problems of the calculus of
variations. First results on Lipschitzian regularity for the
basic problem of the calculus of variations -- $\varphi(u) = u$ --
belong to L.~Tonelli and S.~Bernstein. Some further results
have been obtained by C.~B.~Morrey and more recently by F.~H.~Clarke
and R.~B.~Vinter among others.
For a survey see \cite[Ch.~2]{MR91j:49001} or
\cite[Ch.~11]{MR2001c:49001}. Less is known for the Lagrange problem
of optimal control. Problems whose dynamics is linear and time
invariant -- $\varphi(x,\,u) = A x + B u$ -- were addressed in
\cite{MR91b:49033}. The result is obtained imposing conditions under 
which the problem is reduced into a problem of the calculus of variations
and then using the results available in the literature.
Recently, a new approach to the Lipschitzian regularity has been
developed by A.~Sarychev and the author in \cite{MR2000m:49048},
which allows to deal with a wide class of optimal control problems
with control-affine dynamics -- $\varphi(t,\,x,\,u) = f(t,\,x) + g(t,\,x)\,u$.
As particular cases they include the problems of the calculus of variations 
treated before. The approach is based in the reduction of the problem to a 
time-optimal control problem, on the subsequent compactification of the space of 
admissible controls, and utilization of 
the Pontryagin's maximum principle \cite{MR29:3316b}. 
The conditions of Lipschitzian regularity arise from the conditions
of applicability of the Pontryagin's maximum principle \cite{MR29:3316b} 
to the latter problem
and from its equations Lipschitzian regularity of the corresponding minimizer
of the initial problem is established. For a survey see \cite{MR2001j:49062}.
The main result of \cite{MR2000m:49048} can be summarized in the following theorem:
\begin{theorem}[$\mathbf{\varphi\left(t,x,u\right)= 
f\left(t,x\right)+g\left(t,x\right)\,u}$]
\label{mr:MR2000m:49048}
If $g\left(t,x\right)$ has complete rank $r$
for all $t$ and $x$; $L(\cdot,\cdot,\cdot)$, $f\left(\cdot,\cdot\right)$, 
$g\left(\cdot,\cdot\right)$ are $C^1$ smooth; and the following conditions
are satisfied:
\begin{description}
\item[(coercivity)]  there exist a function 
$\theta:\mathbb{R}\rightarrow\mathbb{R}$ and 
$\zeta\in\mathbb{R}$ such that
\begin{gather*}
L\left(t,x,u\right)\geq \theta\left(\left\|u\right\|\right) >\zeta \, ,
\quad \forall \, \left(t,x,u\right) \, , \\
\lim_{r\rightarrow+\infty}\frac{\theta\left(r\right)}{r}= + \infty \, ;
\end{gather*}
\item[(growth condition)]  there exist constants $\gamma $, $\beta $, 
$\eta $ and $\mu$, with $\gamma >0$, $\beta <2$ and 
$\mu \geq \max \left\{\beta -2,\,-2\right\}$, such that 
the inequality
\begin{equation*}
\left( \left| L_{t}\right| +\left| L_{x^{i}}\right| +\left\| L\,
\varphi_{t}-L_{t}\,\varphi \right\| +\left\| L\,\varphi _{x^{i}}-L_{x^{i}}\,
\varphi\right\| \,\right) \,\left\| u\right\| ^{\mu }\leq \gamma \,L^{\beta }
+\eta \, ,
\end{equation*}
holds for all $t\in \left[a,b\right]$, 
$x\in\mathbb{R}^{n}$, $u\in\mathbb{R}^{r}$,
$i\in \left\{ 1,\,\ldots ,\,n\right\}$;
\end{description}
then all the minimizers $\tilde{u}\left(\cdot\right)$ of the problem, 
which are not abnormal extremal controls, are essentially bounded on 
$\left[a,b\right]$.
\end{theorem}
Results for general nonlinear dynamics,
which is nonlinear both in state and control variables, 
are lacking. To deal with the problem we make use of
a different auxiliary optimal control problem than the
one in \cite{MR2000m:49048} (Section~\ref{S:Sec2}), which is 
obtained using an idea of time reparameterization that proved to be useful in
many different contexts -- see \textrm{e.g.} 
\cite[Sec.~10]{MR28:3353}, \cite{dubmil63,dubmil65},
\cite[Lec.~13]{MR57:3958}, \cite[p.~46]{MR85c:49001},
\cite{MR91c:49060}, 
\cite{MR91f:49032}, \cite{MR93f:49002},
\cite[Ch.~5]{MR99g:49002}, 
\cite[p.~29]{MR2000m:49002}, and \cite{delfimEJC}.
At the core of our proof techniques is the study of the relation
between the minimizers and admissible trajectory-control pairs
of the original and auxiliary problems
(Section~\ref{S:Sec3}) and how the Pontryagin extremals are related 
(Section~\ref{S:Sec4}).
Applying weak necessary conditions than the
Pontryagin's maximum principle \cite{MR29:3316b},
for example the ones found in \cite{MR85m:49002},
to minimizers of the auxiliary problem,
then one is able to obtain, from the established relations between
minimizers and extremals, the desired regularity properties for the 
minimizers of the original problem (Section~\ref{s:mainresult:grr}). 
Examples which possess minimizers according to the existence
theory and to which our results are applicable while previously
known Lipschitzian regularity conditions fail
are provided (Section~\ref{S:Exemplo}).

%%%%%%%%%%%%%%%%%%%%%%%%%%%%%%%%%%%%%%%%%%%%%%%%%

\section{Formulation of Problems $\mathbf{(P)}$,
$\mathbf{(P_{\tau})}$ and $\mathbf{(P_{\tau}[w(\cdot)])}$}
\label{S:Sec2}

We are interested in the study of 
Lipschitzian regularity conditions for
the Lagrange problem of optimal control with
arbitrary boundary conditions. For that is
enough to consider the case when the
boundary conditions are fixed:
$x(a) = A$ and $x(b) = B$.
Indeed, if $\tilde{x}(\cdot)$ is a
minimizing trajectory for a Lagrange problem with any other kind of
boundary conditions, then $\tilde{x}(\cdot)$ is also a minimizing
trajectory for the corresponding fixed boundary problem with
$A = \tilde{x}(a)$ and $B = \tilde{x}(b)$. The
data for our problem is then
\begin{equation}
\label{e:dados}
\left[
\begin{gathered}
a,\,b \in \mathbb{R}\quad (a < b) \\
A,\,B \in \mathbb{R}^n \\
L : \mathbb{R} \times \mathbb{R}^n \times \mathbb{R}^r \longrightarrow \mathbb{R} \\
\varphi : \mathbb{R} \times \mathbb{R}^n \times \mathbb{R}^r \longrightarrow
\mathbb{R}^n \text{.}
\end{gathered}
\right.
\end{equation}
We assume
$L(\cdot,\cdot,\cdot),\,\varphi(\cdot,\cdot,\cdot) \in C$, 
and $\varphi(\cdot,\cdot,u),\,\varphi(\cdot,\cdot,u) \in C^1$.
(Smoothness hypotheses on $L$ and $\varphi$ can be weakened, 
as is discussed later in connection with the Pontryagin 
maximum principle.)
The Lagrange problem of optimal control is defined as follows.\\

\noindent \textbf{Problem $\mathbf{(P)}$.}
\begin{equation*}
I\left[x(\cdot),\,u(\cdot)\right] = \int_a^b 
L\left(t,\,x(t),\,u(t)\right)\,\mathrm{d}t \longrightarrow \min
\end{equation*}
\begin{equation}
\label{e:admp1}
\left[
\begin{gathered}
x(\cdot) \in W_{1,\,1}\left([a,\,b];\,\mathbb{R}^n \right),\,
         u(\cdot) \in L_1\left([a,\,b];\,\mathbb{R}^r\right) \\
\dot{x}(t) = \varphi\left(t,\,x(t),\,u(t)\right), \quad a.e. \  t \in
[a,\,b]\\ x(a) = A,\, x(b) = B\text{.}
\end{gathered}
\right.
\end{equation}\\

\noindent The overdot denotes differentiation
with respect to $t$, while the prime will be used in the sequel 
to denote differentiation with respect to $\tau$. 
To derive conditions assuring that the optimal controls
$\tilde{u}(\cdot)$ of problem $(P)$ are essentially bounded, 
$\tilde{u}(\cdot) \in L_\infty$, two auxiliary problems,
defined with the same data \eqref{e:dados}, will be used.\\

\noindent \textbf{Problem $\mathbf{(P_{\tau})}$.}
\begin{equation*}
J\left[t(\cdot),\,z(\cdot),\,v(\cdot),\,w(\cdot)\right]
= \int_a^b L\left(t(\tau),\,z(\tau),\,w(\tau)\right)\,
v(\tau) \, \mathrm{d}\tau \longrightarrow \min
\end{equation*}
\begin{equation}
\label{e:admp2}
\left[
\begin{gathered}
t(\cdot) \in W_{1,\,\infty}\left([a,\,b];\,\mathbb{R} \right),\,
z(\cdot) \in W_{1,\,1}\left([a,\,b];\,\mathbb{R}^n \right) \\
v(\cdot) \in L_{\infty}\left([a,\,b];\,[0.5,\,1.5] \right),\,
w(\cdot) \in L_1\left([a,\,b];\,\mathbb{R}^r\right)\\
\begin{cases}
t'(\tau) = v(\tau)\\
z'(\tau) = \varphi\left(t(\tau),\,z(\tau),\,w(\tau)\right)\,v(\tau)
\end{cases} \\
t(a) = a,\, t(b) = b\\
z(a) = A,\, z(b) = B \, .\\
\end{gathered}
\right.
\end{equation}

\begin{remark}
The fact that the control variable $v(\cdot)$ takes on its values
in the set $[0.5,\,1.5]$, guarantees that $t(\tau)$ has an
inverse function $\tau(t)$.
\end{remark}

\noindent Problem $(P)$ is actually equivalent to problem $(P_{\tau})$,
in the sense that problem $(P)$ can be formally transformed
into problem $(P_{\tau})$ by considering $t$ as a dependent variable 
and introducing a one to one Lipschitzian transformation 
$[a,b] \ni t \mapsto \tau \in [a,b]$. Both problems have 
the same minimum value and a direct relation between admissible 
state-control pairs (\textrm{cf.} Section \ref{S:Sec3}).

The following problem is the same as problem $(P_{\tau})$ except that
$w(\cdot) \in L_1\left([a,\,b];\,\mathbb{R}^r\right)$ 
is fixed and the functional is
to be minimized only over $t(\cdot)$, $z(\cdot)$ (the state variables)
and $v(\cdot)$ (the control variable).\\

\noindent \textbf{Problem $\mathbf{\left(P_{\tau}[w(\cdot)]\right)}$.}
\begin{equation*}
K\left[t(\cdot),\,z(\cdot),\,v(\cdot)\right]
= \int_a^b F\left(\tau,\,t(\tau),\,z(\tau),\,v(\tau)\right)\,
\mathrm{d}\tau \longrightarrow \min
\end{equation*}
\begin{equation*}
\left[
\begin{gathered}
t(\cdot) \in W_{1,\,\infty}\left([a,\,b];\,\mathbb{R} \right),\,
z(\cdot) \in W_{1,\,1}\left([a,\,b];\,\mathbb{R}^n \right) \\
v(\cdot) \in L_{\infty}\left([a,\,b];\,[0.5,\,1.5] \right)\\
\begin{cases}
t'(\tau) = v(\tau)\\
z'(\tau) = f\left(\tau,\,t(\tau),\,z(\tau),\,v(\tau)\right)
\end{cases} \\
t(a) = a,\, t(b) = b\\
z(a) = A,\, z(b) = B \, ,\\
\end{gathered}
\right.
\end{equation*}
where
$F(\tau,\,t,\,z,\,v) = L\left(t,\,z,\,w(\tau)\right)\,v$,
$f(\tau,\,t,\,z,\,v) = \varphi\left(t,\,z,\,w(\tau)\right)\,v$.

\begin{remark}
Problem $(P_{\tau})$ is autonomous while $(P)$ 
and $\left(P_{\tau}[w(\cdot)]\right)$ are not.
\end{remark}

The relation between problem $(P)$ and problem
$\left(P_{\tau}[w(\cdot)]\right)$ is discussed 
in the following two sections.

%%%%%%%%%%%%%%%%%%%%%%%%%%%%%%%%%%%%%%%%%%%%%%%%%

\section{Relation Between the Solutions of the Problems}
\label{S:Sec3}

Let us begin to determine the relation between admissible
pairs for problem $(P)$ and admissible quadruples for
problem $(P_{\tau})$.

\begin{definition}
The pair $\left(x(\cdot),\,u(\cdot)\right)$ is said to be
\emph{admissible for $(P)$} if all conditions in \eqref{e:admp1}
are satisfied. Similarly, 
$\left(t(\cdot),\,z(\cdot),\,v(\cdot),\,w(\cdot)\right)$ is
said to be \emph{admissible for $(P_{\tau})$} if all conditions in
\eqref{e:admp2} are satisfied.
\end{definition}

Next lemma specifies a passage from an admissible
state-control pair of the problem $(P)$ to an admissible 
state-control pair of the problem $(P_{\tau})$. Importantly,
the values of the two problems are the same.
\begin{lemma}
\label{r:p1impp2} 
Let $\left(x(\cdot),u(\cdot)\right)$ be
admissible for $(P)$. Then, for any function $v(\cdot)$ satisfying
\begin{gather}
v(\cdot) \in L_{\infty}\left([a,\,b];\,
\left[0.5,\,1.5\right]\right)\text{,} \label{e:dttv}\\
\int_a^b v(s)\,\mathrm{d}s = b - a\text{,} \label{e:dtttb}
\end{gather}
$t(\tau) = a + \int_{a}^{\tau}
v(s)\,\mathrm{d}s$,
$z(\tau) = x\left(t(\tau)\right)$ and 
$w(\tau)=u\left(t(\tau)\right)$, 
are such that 
$\left(t(\cdot),z(\cdot),v(\cdot),w(\cdot)\right)$ is admissible for
$(P_{\tau})$. Moreover,
\begin{equation}
\label{e:j1eqj2}
J\left[t(\cdot),\,z(\cdot),\,v(\cdot),\,w(\cdot)\right] =
I\left[x(\cdot),\,u(\cdot)\right]\text{.}
\end{equation}
\end{lemma}

\begin{proof}
All conditions in \eqref{e:admp2} become satisfied:
\begin{itemize}
\item Function $t(\cdot)$ is Lipschitzian:
$\frac{\mathrm{d}t(\cdot)}{\mathrm{d}\tau} = 
v(\cdot) \in L_{\infty}\left([a,\,b];\,[0.5,\,1.5]\right)$;
\item Function $z(\cdot)$ is absolutely continuous
since it is a composition of the absolutely continuous function
$x(\cdot)$ with the strictly monotonous Lipschitzian continuous
function $t(\cdot)$:
\begin{equation}
\label{e:tmono}
\frac{\mathrm{d}t(\tau)}{\mathrm{d}\tau} = v(\tau) > 0\text{;}
\end{equation}
\item Function $w(\cdot)$ is Lebesgue measurable,
$w(\cdot) \in L_1$,
because $u(\cdot)$ is measurable and $t(\cdot)$ is a strictly
monotonous absolutely continuous function;
\item Differentiating $z(\cdot)$ we obtain:
\begin{equation*}
z'(\tau) = \frac{\mathrm{d}z(\tau)}{\mathrm{d}\tau} =
\frac{\mathrm{d}x\left(t(\tau)\right)}{\mathrm{d}t} \, 
\frac{\mathrm{d}t(\tau)}{\mathrm{d}\tau}\text{.}
\end{equation*}
In view of \eqref{e:admp1} and \eqref{e:tmono}, one concludes from
this last equality that
\begin{equation*}
\begin{split}
z'(\tau) & = \varphi\left(t(\tau),\,x(t(\tau)),\,u(t(\tau))\right) \,
v(\tau) \\ 
& =  \varphi\left(t(\tau),\,z(\tau),\,w(\tau)\right)\,v(\tau)\text{;}
\end{split}
\end{equation*}
\item From \eqref{e:dtttb} and from the definition of $t(\tau)$ we have
$t(a) = a$ and $t(b) = b$. It follows that
\begin{gather*}
z(a) = x\left(t(a)\right) = x(a) = A\text{;} \\
z(b) = x\left(t(b)\right) = x(b) = B\text{.}
\end{gather*}
\end{itemize}
It remains to prove equality \eqref{e:j1eqj2}. Since
\begin{align}
\label{e:j2i}
J\left[t(\cdot),\,z(\cdot),\,v(\cdot),\,w(\cdot)\right]
&= \int_{a}^{b} 
L\left(t(\tau),\,z(\tau),\,w(\tau)\right)\,v(\tau)\,\mathrm{d}\tau \notag \\
&= \int_a^b
L\left(t(\tau),\,x(t(\tau)),\,u(t(\tau))\right)\,v(\tau) \,
\mathrm{d}\tau\text{,}
\end{align}
from the change of variable $t(\tau) = t$,
\begin{equation}
\label{e:mudvar}
\left[
\begin{gathered}
\mathrm{d}t = v(\tau) \, \mathrm{d}\tau \\
\tau = a \Leftrightarrow t = a \\
\tau = b \Leftrightarrow t = b\text{,}
\end{gathered}
\right.
\end{equation}
it follows from \eqref{e:j2i} the pretended conclusion:
\begin{equation*}
J\left[t(\cdot),\,z(\cdot),\,v(\cdot),\,w(\cdot)\right]
= \int_a^b L\left(t,\,x(t),\,u(t)\right) \, \mathrm{d}t =
I\left[x(\cdot),\,u(\cdot)\right]\text{.}
\end{equation*}
\end{proof}

The passage established by Lemma~\ref{r:p1impp2} is not unique
because, compared to $(P)$, problem $(P_{\tau})$ has one more
state variable and one more control variable. However, as far
as the right-hand side of the control system for the problem
$(P_{\tau})$ is autonomous, does not depend on $\tau$, the set
of admissible state-control pairs of $(P_{\tau})$ is invariant
under translations of $\tau$. Lemma~\ref{r:p2impp1} asserts that
to each admissible state-control pair for $(P_{\tau})$ there 
corresponds a unique state-control pair for the problem $(P)$ 
with the same value for the cost functionals.
\begin{lemma}
\label{r:p2impp1}
Let $\left(t(\cdot),\,z(\cdot),\,v(\cdot),\,w(\cdot)\right)$
be admissible for $(P_{\tau})$. Then the pair
$\left(x(\cdot),\,u(\cdot)\right) =
\left(z\left(\tau(\cdot)\right),\,w\left(\tau(\cdot)\right)\right)$, 
where $\tau(\cdot)$ 
is the inverse function of $t(\cdot)$, is admissible for $(P)$.
Moreover
\begin{equation}
\label{e:j2igualj1}
I\left[x(\cdot),\,u(\cdot) \right] =
J\left[t(\cdot),\,z(\cdot),\,v(\cdot),\,w(\cdot)\right].
\end{equation}
\end{lemma}

\begin{proof}
Similar arguments to the ones in the proof of Lemma~\ref{r:p1impp2}, 
show that function $x(\cdot)$ is absolutely continuous and
that $u(\cdot) \in L_1\left([a,\,b];\,\mathbb{R}^r\right)$.
Differentiating $x(\cdot)$ we obtain:
\begin{equation*}
\dot{x}(t) = \frac{\mathrm{d}x(t)}{\mathrm{d}t} = 
\frac{\mathrm{d}z\left(\tau(t)\right)}{\mathrm{d}\tau} \, 
\frac{\mathrm{d}\tau(t)}{\mathrm{d}t} \, ,
\end{equation*}
and from \eqref{e:admp2}, and from the fact that 
$\frac{\mathrm{d}\tau(t)}{\mathrm{d}t} =
\frac{1}{v\left(\tau(t)\right)} > 0$,
one concludes that
\begin{equation*}
\dot{x}(t)  = \frac{\varphi\left(t(\tau(t)),\,z(\tau(t)),\,w(\tau(t))\right)\,
v\left(\tau(t)\right)}{v(\tau(t))}
 = \varphi\left(t,\,x(t),\,u(t) \right) \, .
\end{equation*}
As far as $\tau(a) = a$ and $\tau(b) = b$, it comes
$x(a) = z\left(\tau(a)\right) = z(a) = A$,
$x(b) = z\left(\tau(b)\right) = z(b) = B$,
and all conditions in \eqref{e:admp1} are satisfied. 
Equality \eqref{e:j2igualj1} follows by direct calculations 
from the change of variable $\tau(t) = \tau$ and relations
\eqref{e:mudvar}:
\begin{equation*}
\begin{split}
I\left[x(\cdot),\,u(\cdot)\right]  &= \int_{a}^{b} 
L\left(t,\,x(t),\,u(t)\right) \, \mathrm{d}t 
= \int_a^b L\left(t(\tau(t)),\,z(\tau(t)),\,w(\tau(t))\right)\, \mathrm{d}t \\
&= \int_a^b L\left(t(\tau),\,z(\tau),\,w(\tau)\right)\,v(\tau) \, \mathrm{d}\tau 
= J\left[t(\cdot),\,z(\cdot),\,v(\cdot),\,w(\cdot)\right] \, \text{.}
\end{split}
\end{equation*}
\end{proof}

From Lemmas \ref{r:p1impp2} and \ref{r:p2impp1},
the following two corollaries are obvious. They
establish the relation between the minimizers of
problems $(P)$ and $(P_{\tau})$. As a consequence, 
solving the problem $(P)$ turns out to be the same 
as solving the problem $(P_{\tau})$.
\begin{corollary}
\label{r:sol1impsol2}
If $\left(\tilde{x}(\cdot),\,\tilde{u}(\cdot)\right)$
is a minimizer of problem $(P)$, then, for any function
$\tilde{v}(\cdot)$ satisfying \eqref{e:dttv} and \eqref{e:dtttb}
(\textrm{e.g.} $\tilde{v}(\tau) \equiv 1$), the 4-tuple
\begin{equation*}
\left(\tilde{t}(\cdot),\,\tilde{z}(\cdot),\,\tilde{v}(\cdot),\,
\tilde{w}(\cdot)\right)\text{,}
\end{equation*}
defined by
$\tilde{t}(\tau) = a + \int_{a}^{\tau} \tilde{v}(s)\,\mathrm{d}s$,
$\tilde{z}(\tau) = \tilde{x}\left(\tilde{t}(\tau)\right)$,
$\tilde{w}(\tau) = \tilde{u}\left(\tilde{t}(\tau)\right)$,
is a minimizer to problem $(P_{\tau})$.
\end{corollary}

\begin{corollary}
If $\left(\tilde{t}(\cdot),\,\tilde{z}(\cdot),\,\tilde{v}(\cdot),\,
\tilde{w}(\cdot)\right)$
is a minimizer of problem $(P_{\tau})$, then
the pair $\left(\tilde{x}(\cdot),\,\tilde{u}(\cdot)\right)$ 
defined from $\left(\tilde{t}(\cdot),\,\tilde{z}(\cdot),\,\tilde{v}(\cdot),\,
\tilde{w}(\cdot)\right)$ as in Lemma~\ref{r:p2impp1} is a minimizer
to problem $(P)$.
\end{corollary}

Thus, let $\left(\tilde{x}(\cdot),\,\tilde{u}(\cdot)\right)$ 
be a minimizer of problem $(P)$. From Corollary~\ref{r:sol1impsol2}
we know how to construct a minimizer 
$\left(\tilde{t}(\cdot),\,\tilde{z}(\cdot),\,\tilde{v}(\cdot),\,
\tilde{w}(\cdot)\right)$
to problem $(P_{\tau})$. Obviously, 
as far as problem
$\left(P_{\tau}[\tilde{w}(\cdot)]\right)$ is the same
as problem $(P_{\tau})$ except that $\tilde{w}(\cdot)$ is fixed, $\left(\tilde{t}(\cdot),\,\tilde{z}(\cdot),\,\tilde{v}(\cdot)\right)$
furnishes a minimizer to problem 
$\left(P_{\tau}\left[\tilde{w}(\cdot)\right]\right)$.
Choosing $\tilde{v}(\tau) \equiv 1$ we obtain.

\begin{proposition}
\label{p:usesHomogeneity}
If $\left(\tilde{x}(\cdot),\,\tilde{u}(\cdot)\right)$
is a minimizer of problem $(P)$, then the triple
$\left(\tilde{t}(\tau),\,\tilde{z}(\tau),\,\tilde{v}(\tau)\right) = 
(\tau,\,\tilde{x}(\tau),\,1)$ furnishes
a minimizer to problem
$\left(P_{\tau}[\tilde{u}(\cdot)]\right)$.
\end{proposition}

\noindent Proposition~\ref{p:usesHomogeneity} gives a minimizer
to problem $\left(P_{\tau}[\tilde{u}(\cdot)]\right)$,
where $\tilde{u}(\cdot)$ is an optimal control
to the corresponding problem $(P)$. As far as for
the problem $\left(P_{\tau}[\tilde{u}(\cdot)]\right)$
the admissible controls are already bounded,
one can check that the pertinent hypotheses 
on functions $F$ and $f$ (on functions $L$ and $\varphi$)
required by the Pontryagin maximum principle 
are valid, that is, one can obtain conditions under which
$\left(\tilde{t}(\tau),\,\tilde{z}(\tau),\,\tilde{v}(\tau)\right) 
= (\tau,\,\tilde{x}(\tau),\,1)$ is a Pontryagin extremal.
To conclude that the original minimizer 
$\left(\tilde{x}(\cdot),\,\tilde{u}(\cdot)\right)$
of $(P)$ is also an extremal,
one needs to know how the extremals
of the problems are related. This will be addressed
in the next section.

%%%%%%%%%%%%%%%%%%%%%%%%%%%%%%%%%%%%%%%%%%%%%%%%%

\section{Relation Between the Extremals}
\label{S:Sec4}

At the core of optimal control theory is the celebrated
Pontryagin maximum principle. The maximum principle is
a first order necessary optimality condition for
the optimal control problems.
It first appear in the book \cite{MR29:3316b}. Since then,
several versions have been obtained by weakening
the hypotheses. For example, in \cite{MR29:3316b}
it is assumed that functions $L(\cdot,\cdot,\cdot)$ 
and $\varphi(\cdot,\cdot,\cdot)$
are continuous, and have continuous
derivatives with respect to the state variables $x$:
$L(t,\cdot,u),\,\varphi(t,\cdot,u) \in C^1$.
Instead of the continuity assumption of
$L(\cdot,\cdot,\cdot)$ and $\varphi(\cdot,\cdot,\cdot)$,
a version only requiring that  
functions $L(\cdot,x,\cdot)$ and $\varphi(\cdot,x,\cdot)$ are Borel 
measurable can be found in book \cite[Ch.~5]{MR51:8914}.
There, in order to assure the applicability of the maximum principle,
the following assumption is imposed:
there exists an integrable function $\alpha(\cdot)$ defined on $[a,\,b]$
such that the bound
\begin{gather}
\left\|\frac{\partial L}{\partial x}\left(t,x,u(t)\right)\right\| \le
\alpha(t) \label{e:ber:al1}\\
\left\|\frac{\partial \varphi_i}{\partial x}\left(t,x,u(t)\right)\right\| \le
\alpha(t) \label{e:ber:al2}
\end{gather}
($i = 1,\,\ldots,\,n$) 
holds for all $(t,\,x) \in [a,\,b] \times \mathbb{R}^n$.
The existence and integrability of $\alpha(\cdot)$, and
the bound \eqref{e:ber:al1}--\eqref{e:ber:al2},
are guaranteed under the hypotheses that $L$ and $\varphi$
possess derivatives $\frac{\partial L}{\partial x}$ and 
$\frac{\partial \varphi}{\partial x}$ which are continuous
in $(t,x,u)$, and $u(\cdot)$ is essentially bounded
(these are the hypotheses found in \cite{MR29:3316b}).
Alternative hypotheses are the
following growth conditions (see \cite[Sec.~4.4 and p.~212]{MR85m:49002}):
\begin{equation}
\label{e:tonelli:morrey:type}
\left\|\frac{\partial L}{\partial x}\right\| \le
c \left|L\right| + k \, , \quad
\left\|\frac{\partial \varphi_i}{\partial x}\right\| \le
c \left|\varphi_i\right| + k \, ,
\end{equation}
with constants $c$ and $k$, $c > 0$. 
In the context of the Lipschitzian regularity,
conditions~\eqref{e:tonelli:morrey:type} are particularly
important: they are easy to check in practice.
Those who are familiar with
the Lipschitzian regularity conditions for the basic
problem of the calculus of variations ($\varphi(u) = u$), 
will recognize \eqref{e:tonelli:morrey:type} as a generalization 
and a weak version of the
classical Tonelli--Morrey Lipschitzian regularity condition 
(\textrm{cf. e.g.} \cite[\S 7.1]{MR2001j:49062}) where no restriction
is imposed to $\frac{\partial L}{\partial u}$:\\

\noindent \textbf{Tonelli-Morrey Lipschitzian regularity condition.}
If the Lagrangian $L$ satisfies the growth condition
$\left\|\frac{\partial L}{\partial x}\right\| 
+ \left\|\frac{\partial L}{\partial u}\right\| \le
c \left|L\right| + k$, $c > 0$,
then any solution $\tilde{x}(\cdot)$ to the basic problem of the calculus of
variations, in the class of absolutely continuous functions,
is indeed Lipschitzian ($\tilde{u}(\cdot) = \dot{\tilde{x}}(\cdot) \in L_{\infty}$)
and satisfy the Pontryagin maximum principle.\\

From the fact that \eqref{e:tonelli:morrey:type} is a generalization
of the Tonelli--Morrey Lipschitzian regularity condition,
one can guess a link between the applicability conditions of the maximum
principle and the Lipschitzian regularity conditions. The link 
between the applicability conditions of the classical Pontryagin
maximum principle \cite{MR29:3316b} and the Lipschitzian regularity 
conditions for optimal control problems with control-affine dynamics,
was established in \cite{MR2000m:49048} using a reduction of the problem
to a time-minimal control problem. Here, to deal with
general nonlinear dynamics, we use completely different
auxiliary problems and we will need to apply the maximum 
principle under weaker hypotheses than those in \cite{MR29:3316b}.
This is due to the fact that when we fix
$w(\cdot) \in L_1\left([a,\,b];\,\mathbb{R}^r\right)$,
functions $F(\tau,\,t,\,z,\,v)$ and $f(\tau,\,t,\,z,\,v)$
of problem $\left(P_{\tau}\left[w(\cdot)\right]\right)$
are not continuous in $\tau$ but only measurable.
Hypotheses \eqref{e:tonelli:morrey:type} are suitable,
as far as they can be directly verifiable for a given problem.
Weaker hypotheses than \eqref{e:ber:al1}
and \eqref{e:ber:al2} can also be considered. In this respect,
important improvements are obtained 
from the use of nonsmooth analysis.
For example, one can substitute \eqref{e:ber:al1}
and \eqref{e:ber:al2} by the weaker conditions
\begin{gather}
\left|L\left(t,x_1,u(t)\right) - L\left(t,x_2,u(t)\right)\right| \le
\alpha(t) \left\|x_1 - x_2\right\| \label{eq:WC11e12} \\
\left|\varphi_i\left(t,x_1,u(t)\right) - \varphi_i\left(t,x_2,u(t)\right)\right| \le
\alpha(t) \left\|x_1 - x_2\right\| \notag
\end{gather}
and formulate the maximum principle in a nonsmooth setting, in 
terms of generalized gradients (see \cite{MR54:3540,MR85m:49002}).
Proving general versions of the maximum
principle under weak hypotheses is still in progress and 
the interested reader is referred to the recent paper \cite{MR2002e:49040}.

\begin{definition}
\label{d:extremal:P}
Let $\left(x(\cdot),\,u(\cdot)\right)$
be admissible for $(P)$. We say that the quadruple $\left(x(\cdot),\,u(\cdot),\,\psi_0,\,\psi(\cdot)\right)$,
$\psi_0 \in \mathbb{R}^{-}_{0}$ and $\psi(\cdot) \in
W_{1,\,1}\left([a,b];\,\mathbb{R}^n \right)$,
is an \emph{extremal} of $(P)$, if the following two
conditions are satisfied for almost all $t \in [a,\,b]$:
\begin{description}
\item[the adjoint system]
\begin{equation}
\label{eq:sh}
\dot{\psi}(t) = 
- \frac{\partial H}{\partial x}\left(t,\,x(t),\,u(t),\,\psi_0,\,\psi(t)\right)\text{;} 
\end{equation}
\item[the maximality condition]
\begin{equation}
\label{eq:cm}
H\left(t,x(t),u(t),\psi_0,\psi(t)\right)
= \sup_{u \in \mathbb{R}^r} H\left(t,x(t),u,\psi_0,\psi(t)\right)
\text{;}
\end{equation}
\end{description}
where the Hamiltonian equals
\begin{equation*}
H\left(t,\,x,\,u,\,\psi_0,\,\psi\right) = \psi_0\,L(t,\,x,\,u) +
\psi \cdot \varphi(t,\,x,\,u)\text{.}
\end{equation*}
\end{definition}

\begin{definition}
\label{d:extremal:Ptau}
Let $\left(t(\cdot),\,z(\cdot),\,v(\cdot),\,w(\cdot)\right)$
be admissible for $(P_{\tau})$. The 7-tuple $\left(t(\cdot),\,z(\cdot),\,v(\cdot),\,w(\cdot),\,p_0,\,p_t(\cdot),\,p_z(\cdot)\right)$,
$p_0 \in \mathbb{R}^{-}_{0}$, $p_t(\cdot) \in
W_{1,\,\infty}\left([a,b];\,\mathbb{R}\right)$ and $p_z(\cdot) \in
W_{1,\,1}\left([a,b];\,\mathbb{R}^n \right)$,
is said to be an \emph{extremal} of $(P_{\tau})$, if the following two
conditions are satisfied for almost all $\tau \in [a,\,b]$:
\begin{description}
\item[the adjoint system]
\begin{gather}
\label{e:adjsyst:extrPtau}
\begin{cases}
p_{t}'(\tau) = - \dfrac{\partial \mathcal{H}}{\partial 
t}\left(t(\tau),\,z(\tau),\,v(\tau),\,w(\tau),\,p_0,\,p_t(\tau),\,
p_z(\tau)\right)\, \text{,}\\[0.3cm]
p_{z}'(\tau) = - \dfrac{\partial \mathcal{H}}{\partial 
z}\left(t(\tau),\,z(\tau),\,v(\tau),\,w(\tau),\,p_0,\,p_t(\tau),\,p_z(\tau)\right)
\, \text{;}
\end{cases}
\end{gather}
\item[the maximality condition]
\begin{multline}
\label{e:maxcond:extrPtau}
\mathcal{H}\left(t(\tau),\,z(\tau),\,v(\tau),\,w(\tau),\,
p_0,\,p_t(\tau),\,p_z(\tau)\right) \\
= \sup_{\substack{v \in [0.5,\,1.5] \\ w \in \mathbb{R}^r}} \mathcal{H}\left(t(\tau),
\,z(\tau),\,v,\,w,\,p_0,\,p_t(\tau),\,p_z(\tau)\right)\text{;}
\end{multline}
\end{description}
where the Hamiltonian equals
\begin{equation*}
\mathcal{H}\left(t,\,z,\,v,\,w,\,p_0,\,p_t,\,p_z\right) =
\left(p_0\,L\left(t,\,z,\,w\right) + p_t +
p_z \cdot \varphi\left(t,\,z,\,w\right)\right) \,v \text{.}
\end{equation*}
\end{definition}

\begin{remark}
Functions $H$ and $\mathcal{H}$, respectively the Hamiltonians
in Definitions~\ref{d:extremal:P} and \ref{d:extremal:Ptau}, are
related by the following equality:
\begin{equation}
\label{eq:relacaoEntreHamiltonianos}
\mathcal{H}\left(t,\,z,\,v,\,w,\,p_0,\,p_t,\,p_z\right)
= \left(H\left(t,\,z,\,w,\,p_0,\,p_z\right) + p_t\right)\,v \, \text{.}
\end{equation}
From it one concludes that
\begin{gather}
\frac{\partial \mathcal{H}}{\partial t} =
\frac{\partial H}{\partial t}\,v \, \text{,} \label{e:dmcHdt:dHdt} \\
\frac{\partial \mathcal{H}}{\partial z} =
\frac{\partial H}{\partial x}\,v \, \text{.} \label{e:dmcHdz:dHdx}
\end{gather}
\end{remark}

\begin{definition}
An extremal is called \emph{normal} if the cost multiplier
($\psi_0$ in the Definition~\ref{d:extremal:P} and 
$p_0$ in the Definition~\ref{d:extremal:Ptau}) 
is different from zero and \emph{abnormal} if it vanishes. 
\end{definition}

\begin{remark}
As far as the Hamiltonian is homogeneous with respect
to the Hamiltonian multipliers, 
for normal extremals one can always consider, by scaling, that the cost
multiplier takes value $-1$.
\end{remark}

\begin{remark}
The (Pontryagin) maximum principle give
conditions, as those discussed in the introduction of this
section, under which 
to each minimizer of the problem there corresponds
an extremal with Hamiltonian multipliers
not vanishing simultaneously
($\left(\psi_0,\,\psi\right) \neq 0$ 
in the Definition~\ref{d:extremal:P} and 
$\left(p_0,\,p\right) \neq 0$ in the 
Definition~\ref{d:extremal:Ptau}).
\end{remark}

One can expect the set of extremals of problem $(P_{\tau})$
to be richer than the set of extremals of problem $(P)$.
Nevertheless, there is a relationship between the extremals
of the problems. Next lemma shows that to each
extremal of problem $(P)$ there corresponds extremals of problem
$(P_{\tau})$ lying on the zero level of the maximized Hamiltonian
$\mathcal{H}$.

\begin{lemma}
\label{r:prop2ecc}
Let $\left(x(\cdot),\,u(\cdot),\,\psi_0,\,\psi(\cdot)\right)$
be an extremal of $(P)$. 
Then, for any function 
$v(\cdot) \in L_{\infty}\left([a,\,b];\,
\left[0.5,\,1.5\right]\right)$ satisfying
$\int_a^b v(s)\,\mathrm{d}s = b - a$,
the 7-tuple
$\left(t(\cdot),\,z(\cdot),\,v(\cdot),\,
w(\cdot),\,p_0,\,p_t(\cdot),\,p_z(\cdot)\right)$ defined by
\begin{gather*}
t(\tau) = a + \int_{a}^{\tau} v(s)\,\mathrm{d}s\, ,  \\
z(\tau) = x(t(\tau)) \, , \quad w(\tau) = u(t(\tau))\, , \\
p_0 = \psi_0\, , \quad p_z(\tau) = \psi(t(\tau))\, ,  \\
p_t(\tau) =-H\left(t(\tau),\,x(t(\tau)),\,u(t(\tau)),\,
\psi_0,\,\psi(t(\tau))\right) \, 
\end{gather*}
is an extremal of $(P_{\tau})$ with
$\mathcal{H}\left(t(\tau),\,z(\tau),\,v(\tau),\,w(\tau),\,p_0,\,p_t(\tau),
\,p_z(\tau)\right) \equiv 0$.
\end{lemma}

\begin{proof}
From Lemma~\ref{r:p1impp2} we know that such 7-tuple 
is admissible for $(P_{\tau})$. The maximality condition
\eqref{e:maxcond:extrPtau} is trivially satisfied since 
we are in the singular case:
from \eqref{eq:relacaoEntreHamiltonianos}
the Hamiltonian $\mathcal{H}$ vanishes for
$p_t =-H\left(t,\,z,\,w,\,p_0,\,p_z\right)$.
It remains to prove the adjoint system \eqref{e:adjsyst:extrPtau}.
Since $\frac{dH}{dt} = \frac{\partial H}{\partial t}$ along
the extremals (see \textrm{e.g.} \cite{MR29:3316b} or \cite{MR51:8914})
the derivative of $p_t(\tau)$ with respect to $\tau$ is given by
\begin{equation*}
\frac{dp_t}{d\tau} = - \frac{dH}{d\tau} =
- \frac{dH}{dt}\,\frac{dt}{d\tau} =
- \frac{\partial H}{\partial t}\,\frac{dt}{d\tau} =
- \frac{\partial H}{\partial t}\,v\text{.}
\end{equation*}
From relation \eqref{e:dmcHdt:dHdt}
the first of the equalities \eqref{e:adjsyst:extrPtau} is proved:
$p_{t}' = - \frac{\partial \mathcal{H}}{\partial t}$.
Similarly, as far as $p_z(\tau) = \psi\left(t(\tau)\right)$ and
from \eqref{eq:sh}
$\frac{d}{dt}\psi(t) = - \frac{\partial H}{\partial x}$,
it follows from \eqref{e:dmcHdz:dHdx} that
$p_{z}' = \frac{d\psi(t)}{dt}\,\frac{dt}{d\tau} =
- \frac{\partial H}{\partial x}\,v =
- \frac{\partial \mathcal{H}}{\partial z}$.
\end{proof}

It is also possible to construct an
extremal of problem $(P)$ given an extremal of
$(P_{\tau})$.

\begin{lemma}
\label{r:propImpAbn}
Let $\left(t(\cdot),\,z(\cdot),\,v(\cdot),\,
w(\cdot),\,p_0,\,p_t(\cdot),\,p_z(\cdot)\right)$ be an
extremal of $(P_{\tau})$. Then
$\left(x(\cdot),\,u(\cdot),\,\psi_0,\,\psi(\cdot)\right) =
\left(z\left(\tau(\cdot)\right),\,w\left(\tau(\cdot)\right),
\,p_0,\,p_z\left(\tau(\cdot)\right)\right)$ is an extremal
of $(P)$ with $\tau(\cdot)$ the inverse function of $t(\cdot)$.
\end{lemma}

\begin{proof}
From Lemma~\ref{r:p2impp1} we know that the pair
$\left(x(\cdot),\,u(\cdot)\right)$ is admissible
for $(P)$. Direct calculations show that
\begin{equation*}
\dot{\psi} = \frac{\mathrm{d}}{\mathrm{d}t} p_z(\tau)
= \frac{\mathrm{d}p_z(\tau)}{\mathrm{d}\tau} \, 
\frac{\mathrm{d}\tau}{\mathrm{d}t} = 
- \frac{\partial \mathcal{H}}{\partial z} \frac{1}{v} \, \text{.}
\end{equation*}
From \eqref{e:dmcHdz:dHdx} the required adjoint system is
obtained: $\dot{\psi} = - \frac{\partial H}{\partial x}$.
Maximality condition \eqref{e:maxcond:extrPtau} implies that
\begin{multline*}
\mathcal{H}\left(t(\tau),\,z(\tau),\,v(\tau),\,w(\tau),\,
p_0,\,p_t(\tau),\,p_z(\tau)\right) \\
= \sup_{w \in \mathbb{R}^r} \mathcal{H}\left(t(\tau),
\,z(\tau),\,v(\tau),\,w,\,p_0,\,p_t(\tau),\,p_z(\tau)\right)
\end{multline*}
for almost all $\tau \in [a,\,b]$. Given the relation
\eqref{eq:relacaoEntreHamiltonianos} one can write that
\begin{equation*}
H\left(t(\tau),\,z(\tau),\,w(\tau),\,p_0,\,p_z(\tau)\right)
= \sup_{w \in \mathbb{R}^r} 
H\left(t(\tau),\,z(\tau),\,w,\,p_0,\,p_z(\tau)\right)\, \text{.}
\end{equation*}
Putting $\tau = \tau(t)$ we obtain
the maximality condition \eqref{eq:cm}.
\end{proof}

Lemmas~\ref{r:prop2ecc} and \ref{r:propImpAbn}
establish a correspondence between abnormal extremals
of problems $(P)$ and $(P_{\tau})$.

\begin{corollary}
If there are no abnormal extremals of problem $(P)$ then there are no
abnormal extremals of problem $(P_\tau)$. If there are
no abnormal extremals of $(P_\tau)$ then there are also no
abnormal extremals of $(P)$.
\end{corollary}

\begin{definition}
We call a control an \emph{abnormal extremal
control} if it corresponds to an abnormal extremal.
\end{definition}

\begin{proposition}
\label{prop:end:secRelExt}
If $\left(\tilde{x}(\cdot),\,\tilde{u}(\cdot)\right)$
is a minimizer of problem $(P)$ and $\tilde{u}(\cdot)$
is not an abnormal extremal control, then the minimizing
control $\tilde{v} \equiv 1$ of Proposition~\ref{p:usesHomogeneity}
is not an abnormal extremal control too.
\end{proposition}

In the next section we will use Proposition~\ref{prop:end:secRelExt}
to show that, for non-abnormal minimizers,
the Lipschitzian regularity conditions we are looking for,
assuring that the minimizing controls, predicted by Tonelli's 
existence theorem, are indeed bounded, appear from
the applicability conditions of the maximum principle to problem 
$\left(P_{\tau}[\tilde{u}(\cdot)]\right)$.

%%%%%%%%%%%%%%%%%%%%%%%%%%%%%%%%%%%%%%%%%%%%%%%%%

\section{The General Regularity Result}
\label{s:mainresult:grr}

Filippov \cite{MR26:7469} gave the first
general existence theorem for optimal control
(the original paper, in russian, appear in 1959).
There exist now an extensive literature on the existence of
solutions to problems of optimal control. We refer the interested
reader to the book \cite{MR85c:49001} for significant results,
various formulations, and detailed discussions.
Follows a set of conditions, of the type of
Tonelli \cite{JFM45.0615.02}, that guarantee existence of minimizer
for problem $(P)$.\\

\noindent \textbf{``Tonelli'' existence theorem for $\mathbf{(P)}$.}
Problem $(P)$ has a minimizer 
$\left(\tilde{x}(\cdot),\,\tilde{u}(\cdot)\right)$ with
$\tilde{u}(\cdot) \in L_1\left([a,\,b];\,\mathbb{R}^r\right)$,
provided there exists at least one admissible pair,
functions $L(\cdot,\cdot,\cdot)$ and
$\varphi(\cdot,\cdot,\cdot)$ are continuous,
and the following convexity and coercivity conditions hold:
\begin{description}
\item[(convexity)] Functions $L(t,x,\cdot)$ and
$\varphi(t,x,\cdot)$  are convex for all $\left(t,\,x\right)$;
\item[(coercivity)] There exists a function 
$\theta : \mathbb{R}_{0}^{+} \rightarrow \mathbb{R}$, bounded below, such
that
\begin{gather}
L(t,\,x,\,u) \ge \theta\left(\left\|\varphi(t,\,x,\,u)\right\|\right)
\quad \text{for all } (t,x,u) \text{;} \label{e:coercivity1}\\
\lim_{r \rightarrow + \infty} \frac{\theta(r)}{r} = +
\infty\text{;} \label{e:coercivity2} \\
\lim_{\left\|u\right\| \rightarrow + \infty}
\left\|\varphi(t,\,x,\,u)\right\| = + \infty
\quad \text{for all } (t,\,x)\text{.} \label{e:coercivity3}
\end{gather}
\end{description}

\begin{remark}
For the basic problem of the calculus of variations
one has $\varphi = u$ and the theorem above
coincides with the classical Tonelli existence theorem.
\end{remark}

Analyzing the hypotheses of both necessary optimality
conditions and existence theorem, one comes to the conclusion
that the requirements of existence theory do not imply
those of the maximum principle. 
Given a problem, it may happen that the necessary optimality conditions
are valid while existence is not guarantee; or 
it may happen that the minimizers predicted by the existence 
theory fail to be extremals. 
Follows the main results of the paper.

\begin{theorem}
\label{r:mainresult} 
Under the above hypothesis of coercivity, 
all control minimizers $\tilde{u}(\cdot)$ of 
$(P)$, which are not abnormal extremal
controls, are essentially bounded on $[a,\,b]$
if the applicability conditions of the maximum principle
(for example \eqref{e:tonelli:morrey:type}, 
\eqref{e:ber:al1}--\eqref{e:ber:al2} or \eqref{eq:WC11e12})
to functions $F$ and $f$
of problem $\left(P_{\tau}\left[\tilde{u}(\cdot)\right]\right)$
are assured.\footnote{A precise statement of possible regularity 
conditions are found in Theorem~\ref{r:cor1:lipreg}.}
\end{theorem}

\begin{remark}
Convexity is not required in the regularity result of 
Theorem~\ref{r:mainresult} in order to establish
the Lipschitzian regularity of the (non-abnormal)
minimizing trajectories $\tilde{x}(\cdot)$. Convexity 
is only required to establish the existence of minimizers, not
the regularity. This fact is
important since existence theorems without the convexity assumptions
are a question of great interest (see \textrm{e.g.} \cite{MR2000i:49005}
and the references therein).
\end{remark}

Applying the hypotheses \eqref{e:tonelli:morrey:type}
of the maximum principle to functions $F$ and $f$
of problem $\left(P_{\tau}\left[\tilde{u}(\cdot)\right]\right)$, the
following result is trivially obtained.

\begin{theorem}
\label{r:cor1:lipreg}
Under the hypothesis of coercivity, the growth conditions: there exist
constants $c>0$ and $k$ such that
\begin{equation*}
\begin{split}
\left|\frac{\partial L}{\partial t}\right| &\le
c \left|L\right| + k \, , \quad 
\left\|\frac{\partial L}{\partial x}\right\| \le
c \left|L\right| + k \, ,\\
\left\|\frac{\partial \varphi}{\partial t}\right\| &\le
c \left\|\varphi\right\| + k \, , \quad 
\left\|\frac{\partial \varphi_i}{\partial x}\right\| \le
c \left|\varphi_i\right| + k \quad (i = 1,\,\ldots,\,n) \, \text{;}
\end{split}
\end{equation*}
imply that all minimizers $\tilde{u}(\cdot)$ of $(P)$,
which are not abnormal extremal controls,
are essentially bounded on $[a,\,b]$.
\end{theorem}
\noindent In the special cases of a time-invariant linear
control system, or for the problems of the calculus
of variations, Theorem~\ref{r:cor1:lipreg} gives the typical growth
conditions, of the type of Tonelli--Morrey,
obtained in previous results 
\cite{MR86h:49020,MR90k:49006,MR91b:49033}. For the
case of control-affine dynamics, the growth condition
of Theorem~\ref{mr:MR2000m:49048} is less restrictive.

A minimizer $\tilde{u}(\cdot)$ which is not essentially bounded
may fail to satisfy the Pontryagin Maximum Principle. As far
as essentially bounded minimizers are concerned, the Pontryagin 
Maximum Principle is valid.

\begin{theorem}
\label{cor:passouATeorema}
Under the hypotheses of Theorem~\ref{r:cor1:lipreg}, 
all minimizers of $(P)$ are Pontryagin extremals.
Furthermore, it is
valid a weaker form of the Pontryagin maximum
principle in which the adjoint multipliers 
$\psi(\cdot)$ are not required to be absolutely
continuous, but are required instead to be
merely Lipschitzian.
\end{theorem}

\begin{proof} \emph{(Theorem~\ref{r:mainresult})}
Let $\left(\tilde{x}(\cdot),\,\tilde{u}(\cdot)\right)$ be a
minimizer of $(P)$. From Propositions~\ref{p:usesHomogeneity}
and \ref{prop:end:secRelExt} and by the assumptions of the
theorem, we know that there exist absolutely continuous
functions $\tilde{p_t}(\cdot)$ and $\tilde{p_z}(\cdot)$
such that for almost all points $\tau \in [a,\,b]$
\begin{equation*}
v \longmapsto \left[- L\left(\tau,\,\tilde{x}(\tau),\,
\tilde{u}(\tau)\right) + \tilde{p_t}(\tau) +
\tilde{p_z}(\tau) \cdot  
\varphi\left(\tau,\,\tilde{x}(\tau),\,\tilde{u}(\tau)\right)\right] \, v
\end{equation*}
is maximized at $v = 1$ on the interval $[0.5,\,1.5]$. This implies that
\begin{equation}
\label{e:convexfunction:interior}
L\left(\tau,\,\tilde{x}(\tau),\,\tilde{u}(\tau)\right) =
\tilde{p_t}(\tau) + \tilde{p_z}(\tau) \cdot 
\varphi\left(\tau,\,\tilde{x}(\tau),\,\tilde{u}(\tau)\right) \, \text{.}
\end{equation}
Let $\left|\tilde{p_t}(\tau)\right| \le M$ and
$\left\|\tilde{p_z}(\tau)\right\| \le M$ on $[a,\,b]$.
Dividing both sides of inequality \eqref{e:convexfunction:interior}
by $\left\|\varphi\left(\tau,\,\tilde{x}(\tau),\,\tilde{u}(\tau)\right)\right\|$
and using the coercivity hypothesis \eqref{e:coercivity1}, one obtains
\begin{equation*}
\frac{\theta\left(\left\|\varphi\left(\tau,\,\tilde{x}(\tau),\,
\tilde{u}(\tau)\right)\right\|\right)}{\left\|\varphi\left(\tau,\,\tilde{x}(\tau),
\,\tilde{u}(\tau)\right)\right\|} \le M \, 
\frac{1 + 
\left\|\varphi\left(\tau,\,\tilde{x}(\tau),
\,\tilde{u}(\tau)\right)\right\|}{\left\|\varphi\left(\tau,
\,\tilde{x}(\tau),\,\tilde{u}(\tau)\right)\right\|} \, \text{.}
\end{equation*}
The coercivity condition \eqref{e:coercivity2}--\eqref{e:coercivity3}
yields the essential boundedness of $\tilde{u}(\cdot)$ on $[a,\,b]$.
\end{proof}

%%%%%%%%%%%%%%%%%%%%%%%%%%%%%%%%%%%%%%%%%%%%%%%%%

\section{An Example} 
\label{S:Exemplo}

As far as Theorem~\ref{r:cor1:lipreg} is able to cover optimal
control problems with dynamics which is nonlinear both in the
state and in the control variables, plenty of examples 
possessing minimizers according to the existence theory
can be easily constructed for which our result is applicable
while previously known Lipschitzian regularity conditions,
such as those in \cite{MR91b:49033} and \cite{MR2000m:49048}, fail.
Follows one such example with $n = r = 2$.

\begin{example}
\label{exmp:exemplo}
\begin{equation*}
\begin{gathered}
\int_0^1 \left(u_1^{2}(t) + u_2^{2}(t)\right)\,\left(\mathrm{e}^{2\,\left(x_1(t) +
x_2(t)\right)} + 1 \right)\,\mathrm{d}t \longrightarrow \min\\
\left\{
\begin{gathered}
\dot{x_1}(t) = \sqrt{u_1^{2}(t) + u_2^{2}(t)} \\
\dot{x_2}(t) = u_2(t)\,\mathrm{e}^{x_1(t) + x_2(t)}
\end{gathered}
\right.\\
x_1(0) = 0,\,x_1(1) = 1,\,
x_2(0) = 1,\,x_2(1) = 1\text{.}
\end{gathered}
\end{equation*}
\end{example}

Here we have:

\begin{equation*}
\begin{gathered}
L(x_1,x_2,u_1,u_2) = \left(u_1^2 + u_2^2\right)\,\left(\mathrm{e}^{2\,(x_1
+ x_2)} + 1\right)\text{;}\\ \varphi(x_1,x_2,u_1,u_2) =
\begin{bmatrix}
\varphi_1 \\ \varphi_2
\end{bmatrix}
=
\begin{bmatrix}
\sqrt{u_1^2 + u_2^2} \\ u_2\,\mathrm{e}^{x_1 + x_2}
\end{bmatrix}\text{.}
\end{gathered}
\end{equation*}
Clearly, all conditions of Tonelli's existence theorem
are satisfied: an admissible quadruple is
$\left(x_1(t),\,x_2(t),\,u_1(t),\,u_2(t)\right) =
\left(t,\,1,\,1,\,0\right)$;
functions $L(\cdot,\cdot,\cdot,\cdot)$
and $\varphi(\cdot,\cdot,\cdot,\cdot)$ are continuous in
$\mathbb{R}^4$; function $L(x_1,x_2,\cdot,\cdot)$ is strictly convex;
$\varphi(x_1,x_2,\cdot,\cdot)$ is convex;
from the inequality
$L = \left(u_1^2 + u_2^2\right)\,\left(e^{2\,(x_1 + x_2)} +
1\right) \ge u_1^2 + u_2^2 + u_2^2 e^{2\,(x_1 + x_2)}$
we have quadratic coercivity ($\theta(r) = r^2$).
Therefore the problem has a
solution for $x_1(\cdot) \, , x_2(\cdot) \in
W_{1,1}\left([0,\,1];\,\mathbb{R} \right)$ and
$u_1(\cdot)$, $u_2(\cdot) \in L_1\left([0,\,1];\,\mathbb{R}
\right)$. Smooth assumptions on data \eqref{e:dados} are satisfied, since
$L(\cdot,\cdot,\cdot,\cdot)$ and $\varphi(\cdot,\cdot,u_1,u_2)$ are of
class $C^{\infty}$. Theorem~\ref{r:cor1:lipreg} allow us to conclude 
that all minimizing controls,
which are not abnormal extremal controls, are bounded:
\begin{itemize}
\item The conditions on $\frac{\partial L}{\partial t}$
and $\frac{\partial \varphi}{\partial t}$ are
trivially satisfied as far as the problem is autonomous:
$L$ and $\varphi$ do not depend explicitly on the time
variable.
\item The growth conditions on 
$\frac{\partial L}{\partial x}$
and $\frac{\partial \varphi}{\partial x}$ are
also satisfied:
$\frac{\partial L}{\partial x_1} = 
\frac{\partial L}{\partial x_2} =
2 e^{2(x_1 + x_2)} \left(u_1^2 + u_2^2\right) \le 2 L$;
$\frac{\partial \varphi_1}{\partial x_1} = 
\frac{\partial \varphi_1}{\partial x_2} = 0$;
$\frac{\partial \varphi_2}{\partial x_1} = 
\frac{\partial \varphi_2}{\partial x_2} = \varphi_2$.
\end{itemize}
Unbounded minimizers, if there are any, are abnormal
extremal controls. From Theorem~\ref{cor:passouATeorema} 
all minimizing controls of the problem (normal or abnormal)
can be identified via the Pontryagin maximum principle.

%%%%%%%%%%%%%%%%%%%%%%%%%%%%%%%%%%%%%%%%%%%%%%%%%

\section{Final Remarks}

In this paper we study properties of minimizing
trajectories for general problems of optimal control
in the cases where controls are unconstrained (like
in the calculus of variations).
We provide conditions which guarantee Lipschitzian regularity
of the minimizing trajectories for the Lagrange problem of
optimal control in the general nonlinear case. These conditions
solve the discrepancy between the optimality and
existence results, assuring that minimizers predicted
by the existence theory satisfy the optimality conditions.
At the same time, undesirable phenomena, like the Lavrentiev
one, are naturally precluded. We show that the conditions
of Lipschitzian regularity are related with the applicability
conditions of Pontryagin's maximum principle. To deal with dynamics
which are control-affine, 
the classical Pontryagin maximum principle \cite{MR29:3316b} is enough
(see \cite{MR2000m:49048}). To treat the general case, a maximum principle under weak
assumptions, like the one in \cite{MR51:8914}, is necessary.
Our approach is based on the
relationship of the extremals of the Lagrange problem with the extremals of 
an auxiliary problem, and on the subsequent utilization of Pontryagin's maximum
principle to the later problem. The maximality condition of
Pontryagin's maximum principle together with the coercivity
assumption of the existence theorem imply the Lipschitzian regularity
of the corresponding minimizer of the original problem.
This approach allows us to deal with more general class
of problems of optimal control with nonlinear dynamics.

It remains to clarify the interconnection between Lipschitzian
regularity and abnormal extremality. For the problems of the calculus
of variations studied in \cite{MR86h:49020} and \cite{MR90k:49006}
no abnormal extremals exist. For the optimal control problems considered
in \cite{MR91b:49033} and \cite{MR2000m:49048}, abnormal extremals are, like
here, put aside. The question of how
to establish Lipschitzian regularity for the abnormal
minimizing trajectories seems to be a completely open question.

%%%%%%%%%%%%%%%%%%%%%%%%%%%%%%%%%%%%%%%%%%%%%%%%%

\section*{Acknowledgements}

This work is part of the author's Ph.D. project 
\cite{torresPhD}, carried
out at the University of Aveiro, Portugal, under scientific 
supervision of A.~V.~Sarychev. The author is grateful to him
for the inspiring conversations and for the many helpful comments
on the present topic.

%%%%%%%%%%%%%%%%%%%%%%%%%%%%%%%%%%%%%%%%%%%%%%%%%

\bibliographystyle{abbrv}
\bibliography{torres}

%%%%%%%%%%%%%%%%%%%%%%%%%%%%%%%%%%%%%%%%%%%%%%%%%

\end{document}